\documentclass{ifacconf}

\usepackage{graphicx,algorithm,algorithmic,amsfonts}      
\usepackage{natbib}        
\usepackage{srcltx}

\usepackage{amssymb, amsmath, epsfig}
\usepackage{color}
\usepackage{booktabs}

\usepackage{epstopdf}       

\graphicspath{{./}{figs/}}

\newcommand{\R}{\mathbb R}

\renewcommand{\S}{\mathcal S\,}

\begin{document}
\begin{frontmatter}

\title{A residual based snapshot location strategy for POD in distributed optimal control \\of linear parabolic equations}
\thanks[footnoteinfo]{We would like to thank Z J. Zhou, from Shandong Normal University, China, for providing the data and code of the space-time approximation in \cite{GHZ12}.}

\author[First]{A. Alla} 
\author[First]{C. Gr\"a\ss{}le}
\author[First]{M. Hinze} 

\address[First]{Universit\"at Hamburg, 
  Bundesstr. 55, 20146, Hamburg, Germany \\(e-mail: \{alessandro.alla, carmen.graessle, michael.hinze\}@uni-hamburg.de).}


\begin{abstract}
In this paper we study the approximation of a distributed optimal control problem for linear para\-bolic PDEs with model order reduction based on Proper Orthogonal Decomposition (POD-MOR). POD-MOR is a Galerkin approach where the basis functions are obtained upon information contained in time snapshots of the parabolic PDE related to given input data. In the present work we show that for POD-MOR in optimal control of parabolic equations it is important to have knowledge about the controlled system at the right time instances. For the determination of the time instances (snapshot locations) we propose an a-posteriori error control concept which is based on a reformulation of the optimality system of the underlying optimal control problem as a second order in time and fourth order in space elliptic system which is approximated by a space-time finite element method.  
Finally, we present numerical tests to illustrate our approach and to show the effectiveness of the method in comparison to existing approaches.
\end{abstract}


\end{frontmatter}

\section{Introduction}

Optimization with PDE constraints is nowadays a well-studied topic motivated by its relevance in industrial applications. We are interested in the numerical approximation of such optimization problems in an efficient and reliable way using surrogate models obtained with POD-MOR. The surrogate models in the present work are in general built upon snapshot information of the system. This idea was introduced in \cite{Sir87}.\\
\noindent
Several works focus their attention on the choice of the snapshots in order to approximate either dynamical systems or optimal control problems by suitable surrogate models. Here we mention the work related to the computation of the snapshots for dynamical system in  \cite{KV10} and \cite{HL14}.\\
In optimal control problems the reduced model is usually built upon a forecast on the control which in general does not guarantee a proper construction of the surrogate model, since one does not know how far away the optimal solution is from this forecast control. More sophisticated approaches select snapshots by solving an optimization problem in order to improve the selection of the snapshots according to the desired controlled dynamics. For this purpose optimality system POD (shortly OS-POD) is introduced in \cite{KV08}. Adaptive adjustments of the surrogate models are proposed in \cite{AH01} and in \cite{AFS02}.
In our paper, we address the question of efficient selection of snapshot locations by means of an a-posteriori error control approach proposed in \cite{GHZ12}, where the optimality system is rewritten as a second order in time and a fourth order in space elliptic equation. In particular, the time-adaptivity is used to build the snapshot grid which should be used to construct the POD-MOR surrogate model for the approxi\-mate solution of the optimal control problem. Here the contribution for the reduced control problem is twofold: we directly obtain snapshots related to an approximation of the optimal control and, at the same time, obtain information about the time grid. \\
\noindent
The outline of this paper is as follows. In Section 2 we present the optimal control problem together with the optimality conditions. In Section 3 we present the main results of \cite{GHZ12}. POD and its application to optimal control problems is presented in Section 4. Here, we explain our strategy for snapshot location. Finally, numerical tests are discussed in Section 5. A conclusion and outlook are given in Section 6.

\section{Optimal Control Problem}

In this section we describe the distributed optimal control problem.  The governing equation is given by

\begin{equation}\label{heat}
\left.
\begin{array}{rcll}
y_t-\nu\Delta y & = & f+u &\text{ in } \Omega_T,\\
y(\cdot,0) & = & y_0 & \text{ in } \Omega,\\
y &= & 0 &\text{ on } \Sigma_T,
\end{array}
\right\}
\end{equation}
\noindent
where $\nu, T$ are given real positive constants, $\Omega_T:=\Omega\times (0,T]$, $\Omega\subset\R^n$ is an open bounded domain with smooth boundary, $\Sigma_T:=\partial\Omega\times (0,T]$. The state and the control are denoted by  $ y:\Omega_T\rightarrow\R$ and  $u:\Omega_T\rightarrow\R$, respectively.
We note that \eqref{heat} for $y_0\in L^2(\Omega)$ and $u\in L^2(0,T;H^{-1}(\Omega))$ admits a unique solution $y=y(u)\in W(0,T)$, where
$$W(0,T):=\left\{v\in L^2\left(0,T;H^1_0(\Omega)\right), \dfrac{\partial v}{\partial t}\in L^2\left(0,T;H^{-1}(\Omega)\right)\right\}.$$

\noindent
The cost functional we want to minimize is given by
\begin{equation}
J(y,u):=\dfrac{1}{2} \|y-y_d\|^2_{L^2(\Omega_T)}+\dfrac{\alpha}{2} \|u\|^2_{L^2(\Omega_T)},
\end{equation}
where $y_d\in L^2(\Omega_T)$ is the desired state and $\alpha$ is a real positive constant.
The optimal control problem, then, can be formulated as
\begin{equation}\label{ocp}
\min_{u\in L^2(\Omega_T)} J(y(u),u) \mbox{  where } y(u) \mbox{ satisfies } \eqref{heat}.
\end{equation}
It is easy to argue that \eqref{ocp} admits a unique solution $u\in L^2(\Omega_T)$ with associated state $y(u)\in W(0,T)$, see e.g. \cite{Tro10}.

\noindent
The optimality system of the optimal control problem \eqref{ocp} is given by the state equation \eqref{heat} together with the adjoint equation
\begin{equation}\label{adj}
\left.
\begin{array}{rcll}
-p_t-\nu\Delta p & = & y-y_d &\text{ in } \Omega_T,\\
p(\cdot,T) & = & 0 &\text{ in } \Omega,\\
p & = & 0 &\text{ on } \Sigma_T,
\end{array}
\right\}
\end{equation}
and the optimality condition
\begin{equation}\label{opt_con}
\alpha u + p=0\quad \mbox{ in }\Omega_T.
\end{equation}
\noindent
Since our domain is smooth, the regularities of the optimal state, the optimal control and the associated adjoint state are limited through the regularities of the initial state $y_0$ and of the desired state $y_d.$\\
\noindent
The numerical approximation of the optimality system \eqref{heat}-\eqref{adj}-\eqref{opt_con} with a standard Finite Element Method (FEM) leads to a high-dimensional boundary value problem with respect to time:

\begin{equation}\label{opt_disc}
\left.
\begin{array}{rclrcl}
M\dot{y}^N-\nu A y^N & = & f^N+ u^N, & \quad y^N(0) & = & y_0^N.\\
-M\dot{p}^N-\nu A p^N & = & y^N-y^N_d, & \quad p^N(T) & = & 0.\\
 \alpha u^N + p^N & = & 0 & & & 
\end{array}
\right\}
\end{equation}
Here $y^N, p^N:[0,T]\rightarrow\R^N$ are the 
semi-discrete state and adjoint, respectively,  $\dot{y}^N$ is the time derivative,
$M\in\R^{N\times N}$ denotes the finite element mass matrix and 
$A\in\R^{N\times N}$ the finite element stiffness matrix. Note that the 
dimension $N$ of each equation in the semi-discrete system \eqref{opt_disc} is related to the number of element nodes chosen in the FEM approach.


\section{Space-Time approximation}

In this section we transform the first order optimality conditions \eqref{heat}-\eqref{adj}-\eqref{opt_con} into an elliptic equation of fourth order in space and second order in time involving only the state variable $y$. 
The interested reader can find more details on this section in \cite{GHZ12}. Here, we recall the main results. \\
Under natural assumptions the optimality system \eqref{heat}-\eqref{adj}-\eqref{opt_con} may be rewritten as a boundary value problem in space-time:
\begin{equation}\label{2ord}
\left.
\begin{array}{rcll}
-y_{tt}+\nu\Delta^2 y+\dfrac{1}{\alpha}y & = &\dfrac{1}{\alpha}y_d & \text{ in } \Omega_T,\\
y(\cdot,0) & = & y_0 &\text{ in } \Omega,\\
y & = & 0 &\text{ on } \Sigma_T,\\
\Delta y & = & 0 &\text{ on } \Sigma_T,\\
\left(y_t-\Delta y\right)(T) & = & 0 &\text{ in }\Omega.
\end{array}
\right\}
\end{equation}
Here we set $f\equiv0$. In \cite{GHZ12}, under suitable assumptions on the data, it is shown that the solution to problem \eqref{ocp} 
satisfies \eqref{2ord} almost everywhere in space time.\\
Next we provide the weak formulation of \eqref{2ord}. 
For this purpose let:
$$H^{2,1}_0(\Omega_T):=\left\{v\in H^{2,1}(\Omega_T): v(0)=0 \mbox { in }\Omega\right\},$$
where 
$$H^{2,1}(\Omega_T):=L^2\left(0,T;H^2(\Omega\right)\cap H^1_0\left(\Omega)\right)\cap H^1\left(0,T; L^2(\Omega)\right),$$ and is equipped with the norm
$$\|w\|_{H^{2,1}(\Omega_T)}^2:=\left(\|w\|^2_{L^2(0,T;H^2(\Omega))}+\|w\|^2_{H^1(0,T;L^2(\Omega))}\right).$$
\noindent
The bilinear form and the linear form in the weak formulation are defined as:
$$A_T:H^{2,1}_0(\Omega_T)\times H^{2,1}_0(\Omega_T)\rightarrow\R,\qquad L_T:H^{2,1}_0(\Omega_T)\rightarrow\R,$$
$\qquad A_T(v,w):=\displaystyle\int_{\Omega_T} \left( v_tw_t+\dfrac{1}{\alpha}vw\right)$\\
\hspace*{3cm}$+\displaystyle\int_{\Omega_T}\Delta v\Delta w +\displaystyle\int_\Omega \nabla v(T) \nabla w(T)$,\\
$\qquad L_T(v):=\displaystyle\int_{\Omega_T}\dfrac{1}{\alpha}y_dv.$\\
The weak formulation of equation \eqref{2ord}  then reads: 
given $y_d\in L^2(\Omega_T),\, y_0\in H^1_0(\Omega),$ find $ y\in H^{2,1}(\Omega_T)$ with $y(0)=y_0$ and
\begin{equation}\label{weak_2ord}
A_T(y,v)=L_T(v)\quad \forall v\in H^{2,1}_0(\Omega_T).
\end{equation}
Existence and uniqueness of the solution is proved e.g. in \cite{NPS09} and \cite{GHZ12}.

\noindent

We put our attention on the discrete approximation of \eqref{2ord}. As a first step 
we recall a-priori and a-posteriori error estimates for the semi-discrete problem, 
where the space is kept continuous. Let us consider the time 
discretization $0=t_0<t_1<\ldots<t_n=T$ with $\Delta t_j=t_j-t_{j-1}$ 
and $\Delta t:=\max_j \Delta t_j$. Let $I_j:=[t_{j-1},t_j]$, we define
$$V_t^k:=\left\{v\in H^{2,1}(\Omega_T): \; v(x,\cdot)|_{I_j} \in P_1(I_j) \text{ for a.e. } x \in \Omega \right\},$$
and $$ \bar{V}_t^k:=V_t^k\cap  H^{2,1}_0(\Omega_T).$$
We now consider the semi-discrete problem
\begin{equation}\label{weak_dis}
A_T(y_k,v_k)=L_T(v_k),\quad \forall v_k\in \bar{V}_t^k.
\end{equation}
Under suitable assumptions on the data we have the following a-priori error estimate from \cite[Thm 3.1]{GHZ12}:
$$\|y-y_k\|_{L^2(\Omega_T)} + \Delta t \|y-y_k\|_{H^1(0,T;L^2(\Omega))} \leq C \Delta t^2.$$
Furthermore, from \cite[Thm 3.3]{GHZ12} we also have the temporal a posteriori error 
estimate
\begin{equation}\label{est-thm31}
\|y-y_k\|_{H^{2,1}_0(\Omega_T)}^2\leq C\eta^2,
\end{equation}
where 
\begin{eqnarray}\label{res:2ord}
\eta^2&=&\sum_j \Delta t_j^2 \int_{I_j} \| \dfrac{1}{\alpha} y_d+\dfrac{\partial^2 y_k}{\partial t^2}-\dfrac{1}{\alpha}y_k-\Delta^2 y_k\|_{L^2(\Omega)}^2 \nonumber\\
&+& \sum_j \int_{I_j} \|\Delta y_k\|_{L^2(\Gamma)}^2.
\end{eqnarray}

With the help of \eqref{est-thm31}, we are able to refine the time grid by means of the residual of the system \eqref{2ord}. The snapshot selection for POD-MOR of problem \eqref{ocp} will now be based on \eqref{est-thm31}.


\section{POD for optimal control problems}

In this section we explain the POD method which we utilize in order to replace 
the original (high-dimensional) problem \eqref{ocp} by a (low-dimensional) 
POD-Galerkin approximation. The main interest when applying the POD method is to reduce calculation times and storage capacity while retaining a 
satisfactory approximation quality. This is made possible due to the key fact that POD basis functions contain information about the 
underlying model, since the POD basis are derived from snapshots of 
a solution data set. For this reason it is important to use rich snapshot ensembles reflecting 
the dynamics of the modelled system. 
The snapshot form of POD proposed by \cite{Sir87} works in the discrete version as follows.\\
\noindent
Let us suppose that the continuous solution $y(t)$ of \eqref{heat} belongs to a real separable Hilbert space $V$, equipped with the inner product $\langle\cdot,\cdot\rangle$. In fact, in the applications we use either $V=L^2(\Omega)$ or $V=H^1_0(\Omega)$. For 
given time instances $t_0, t_1, \ldots, t_n$ we suppose to know the solution $y(t_j)$ of \eqref{heat} for $j=0, \ldots, n$. We define the 
snapshot set $\mathcal{V}:=\mbox{span}\{y(t_0),\ldots, y(t_n)\} \subset V$ and determine a POD basis $\{\psi_1,\ldots,\psi_\ell\}$ of 
rank $\ell \in \{1, ..., d\}, d= \dim \mathcal{V} \leq N$, by solving the following minimization problem:
\begin{eqnarray}\label{prb:pod}
\min_{\psi_1,\ldots,\psi_\ell \in \mathbb{R}^N} \sum_{j=0}^n \beta_j \left\|y(t_j)-\sum_{i=1}^\ell \langle y(t_j),\psi_i\rangle \; \psi_i\right\|^2\nonumber\\
\mbox{ s.t. } \langle\psi_j,\psi_i\rangle=\delta_{ij}\quad\mbox{for }1\leq i,j\leq \ell,
\end{eqnarray}
where $\beta_j$ are nonnegative weights and $\delta_{ij}$ denotes the Kronecker symbol. The associated norm is given by \mbox{$\|\cdot\|^2=\langle\cdot,\cdot\rangle$.} \\
\noindent
It is well-known (see e.g. \cite{GV13}) that a solution to problem 
(\ref{prb:pod}) is given by the first $\ell$ eigenvectors $\{\psi_1, \ldots, 
\psi_\ell\}$ corresponding to the $\ell$ largest eigenvalues
$\lambda_i$ of the self-adjoint linear operator $\mathcal R: V \rightarrow V,$ i.e. $\mathcal{R}\psi_i=\lambda_i\psi_i$ with $\lambda_i>0$, where $\mathcal R$ is defined as follows:
$$ \mathcal{R}\psi=\sum_{j=0}^n \beta_j \langle y(t_j),\psi\rangle \; y(t_j) \quad \mbox{for } \psi\in  V.$$ 
\noindent
Moreover, we can quantify the POD approximation error by the neglected eigenvalues depending on the snapshots (more details can be found in e.g. \cite{GV13}) as follows:
\begin{equation}\label{err-POD}
 \sum_{j=0}^n \beta_j \left\|y(t_j)-\sum_{i=1}^\ell \langle y(t_j),\psi_i\rangle \; \psi_i\right\|^2=\sum_{i=\ell+1}^d \lambda_i.
 \end{equation}

\noindent
Let us assume we are able to compute POD basis functions 
for the optimal control problem \eqref{ocp}. Then, we 
define the POD ansatz of order $\ell$ for the state $y$ as
\begin{equation}\label{pod_ans}
y^\ell(t)=\sum_{i=1}^\ell w_i(t)\psi_i,
\end{equation}

\noindent
where $y^\ell\in V^\ell\subset V$, the POD basis functions are 
$\{\psi_i\}_{i=1}^\ell$ and the unknown coefficients are $\{w_i\}_{i=1}^\ell$.
If we plug this assumption into \eqref{heat} we get the following reduced system
 \begin{equation}\label{opt_pod:1} 
 M^\ell\dot{w}-\nu A^\ell w  = u^\ell,  \quad w(0) = w_0,
 \end{equation}
where entries of the mass matrix $M^\ell$ and the stiffness matrix $A^\ell$ are given by $(M^\ell)_{ij}=\langle\psi_j, 
 \psi_i\rangle$ and $(A^\ell)_{ij} =  \langle\nabla\psi_i,\nabla\psi_j\rangle$,  respectively.  
 The coefficients of the initial condition 
 $y^\ell(0)\in\R^\ell$ are determined by $w_i(0)=(w_0)_i=
 \langle y_0,\psi_i\rangle, \;\;1\leq i\leq \ell,$ and the 
projected desired state is obtained as 
 $(w_d)_i=\langle y_d,\psi_i\rangle, 1 \leq i \leq \ell$.

%
Thus, the optimal control problem of reduced 
order is obtained through replacing \eqref{ocp} by a dynamical system obtained from a 
Galerkin approximation with modes 
$\{\psi_i\}_{i=1}^\ell$ and ansatz \eqref{pod_ans} for the 
state variable. This POD-Galerkin approximation leads 
to  the optimization problem
\begin{equation}\label{ocp_pod}
\min_{u \in L^2(\Omega_T)} \hat{J}^\ell(u) \mbox{ s.t. } y^\ell(u) \mbox{ satisfies } \eqref{opt_pod:1},
\end{equation}
where $\hat{J}^\ell$ is the reduced cost functional, i.e. $\hat{J}^\ell 
(u):= \hat{J} (y^\ell(u),u) $. We recall that the discretization of the optimal solution $\bar{u}^\ell$ to \eqref{ocp_pod} is determined 
by the relation between the adjoint state and control and refer 
to \cite{H05} for more details about this variational discretization concept. \\

 Now let us draw our attention to our {\em snapshot 
 location} strategy for POD model order reduction in optimal control. 
 Recalling the minimization problem \eqref{prb:pod}, we observe the strong dependence of the POD basis functions 
 on the chosen snapshots. These snapshots shall have the 
 property to capture the main features of the dynamics of the truth solution
 as good as possible. Here, we face some difficulties since the reduction of optimal control 
 problems is usually initialized with snapshots computed from a given input control $u_{\text{sg}}$ ('$u_{\text{sg}}$': control 
 input for \textbf{s}napshot \textbf{g}eneration). In general, we 
 do not have any information about the optimal control, so that 
in POD-MOR often $u_{\text{sg}}\equiv 0$ is chosen.\\
The a-posteriori error control concept for \eqref{2ord} now offers the possibility to select snapshot locations by a time-adaptive procedure. For this purpose, \eqref{2ord} is solved adaptively in time, where the spatial resolution ($\Delta x$ in Algorithm \ref{Alg:OPTPOD}) is chosen very coarse. The resulting time grid points now serve as snapshot locations, on which our POD-MOR model for the optimization is based, where the snapshots are now obtained from a simulation of \eqref{heat} with high spatial resolution $h$. The right-hand side $u$ is obtained from \eqref{opt_con} with $p$ from \eqref{adj} computed with spatially coarse resolution. The procedure is summarized in Algorithm \ref{Alg:OPTPOD}.

\begin{algorithm}[htbp]
\caption{Adaptive snapshot selection for optimal control problems.}
\label{Alg:OPTPOD}
\begin{algorithmic}[1]
\REQUIRE coarse spatial grid size $\Delta x$, fine spatial grid size $h$, max number of degrees of freedom (dof) for the adaptive time discretization, $T>0$.
\STATE Solve \eqref{2ord} adaptively w.r.t. time with spatial resolution $\Delta x$ and obtain time grid $\mathcal{T}$ with solution $y_{\Delta x}$. 
\STATE Solve \eqref{adj} on $\mathcal{T}$ with spatial resolution $\Delta x$ and r.h.s. $y_{\Delta x}$ to obtain $p_{\Delta x}$. Set $u=-\dfrac{1}{\alpha}p_{\Delta x}.$
\STATE Solve \eqref{heat} on $\mathcal{T}$ with spatial resolution $h$ and r.h.s. $u$ extended to the spatial grid with resolution $h$.
\STATE Compute a POD basis of order $\ell$ and build the POD-MOR model.
\STATE Solve \eqref{ocp_pod} with time grid $\mathcal{T}$ and POD-MOR model \eqref{opt_pod:1}.
%
%
\end{algorithmic}
\end{algorithm}

\section{Numerical Test}

In our numerical computations we use a one-dimensional spatial 
domain and a finite element discretization in space by means of conformal 
piecewise linear polynomials with an implicit Euler discretization in time.
The coefficients $\beta_j$ in \eqref{prb:pod} are chosen according to the trapezoidal rule.
All coding is done in \textsc{Matlab R2015}a, and the 
computations are performed on a 2.50GHz computer.\\
\noindent
In the following numerical test, we apply the adaptive snapshot selection 
strategy of Algorithm \ref{Alg:OPTPOD} to demonstrate which influence the chosen time grid can 
have on the approximation quality of the POD suboptimal 
solution and compare with results obtained on an equidistant time grid. 

\textbf{Test Example 5.1.} The data for 
this test example is taken from Example 5.2 in \cite{GHZ12}, where the 
following choices are made: $\Omega = (0,1)$ and 
$[0,T] = [0,1]$. We set $\alpha = 1$ and $\nu = 1$. The example is 
built in such a way that the exact optimal 
solution $(\bar{y},\bar{u})$ 
of problem \eqref{ocp} with associated optimal 
adjoint state $\bar{p}$ is 
known: $$\bar{y}(x,t) = \sin (\pi x) \text{atan} \left(\dfrac{t-1/2}{\varepsilon}\right),$$ \;
$$\bar{u}(x,t) = -\sin (\pi x) \sin (\pi t), \; \bar{p}(x,t) = \sin (\pi x) \sin (\pi t).$$ 

\begin{figure}[htbp]
 \centering
\includegraphics[scale=0.5]{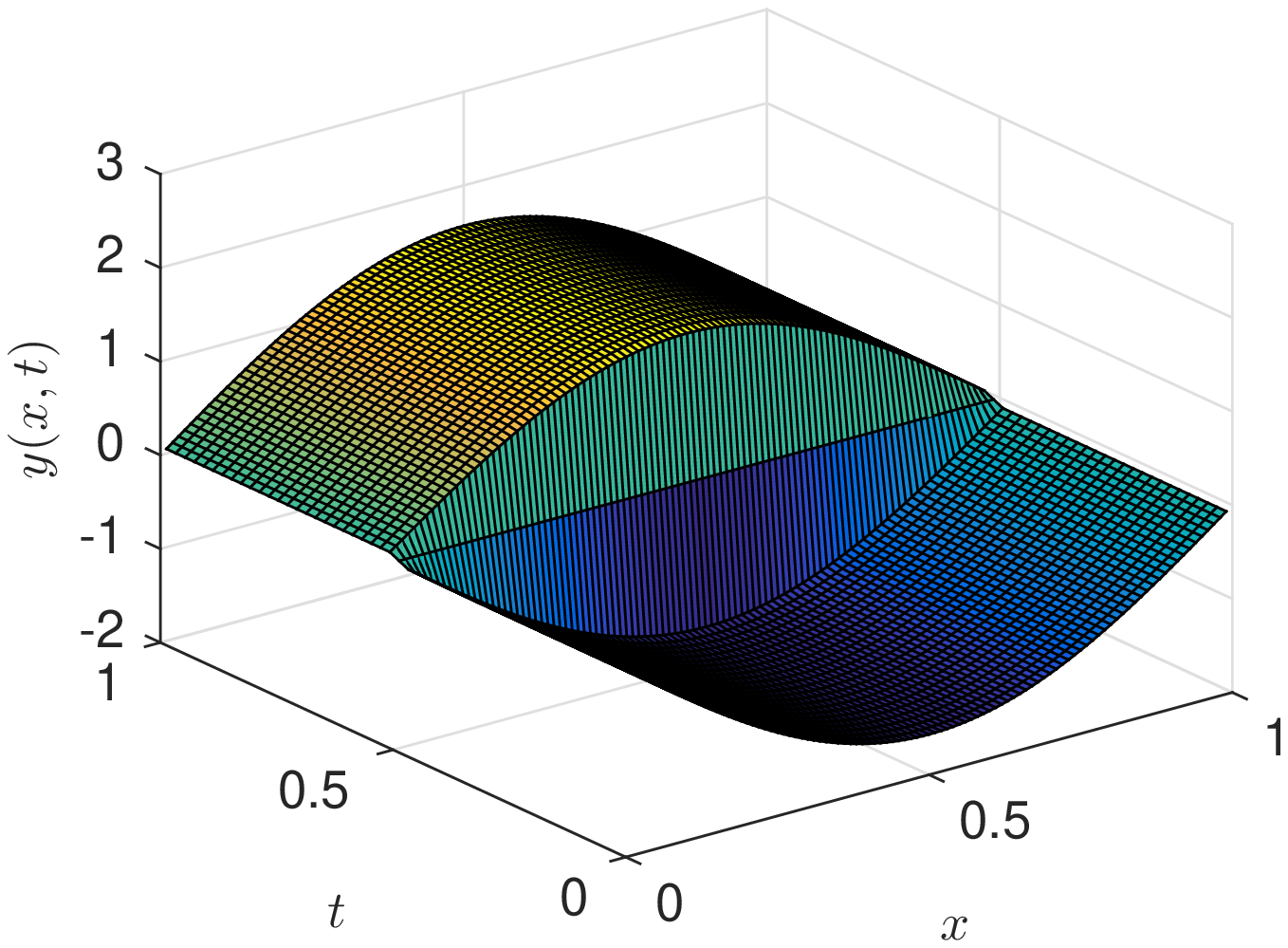} \includegraphics[scale=0.5]{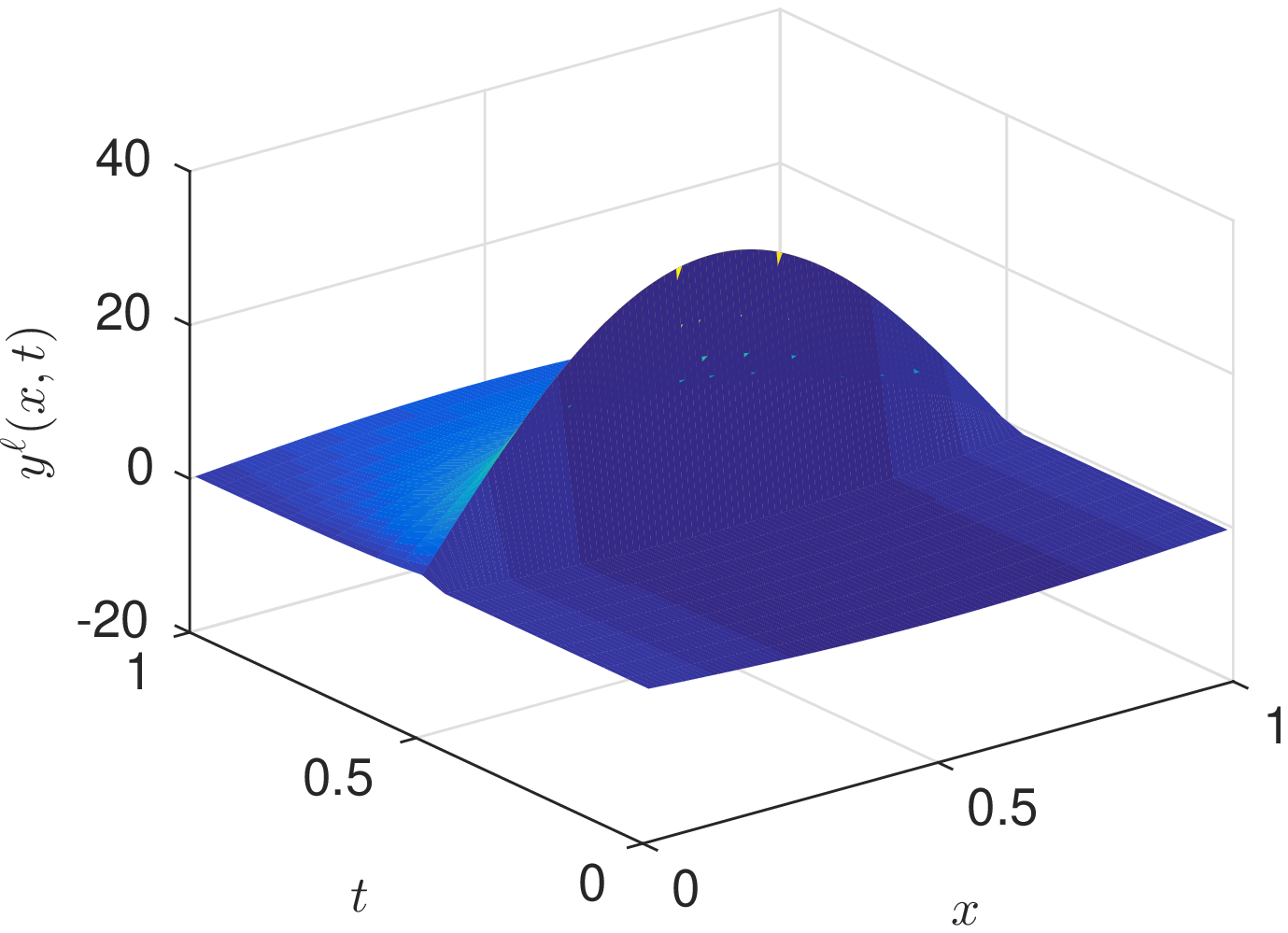} \includegraphics[scale=0.5]{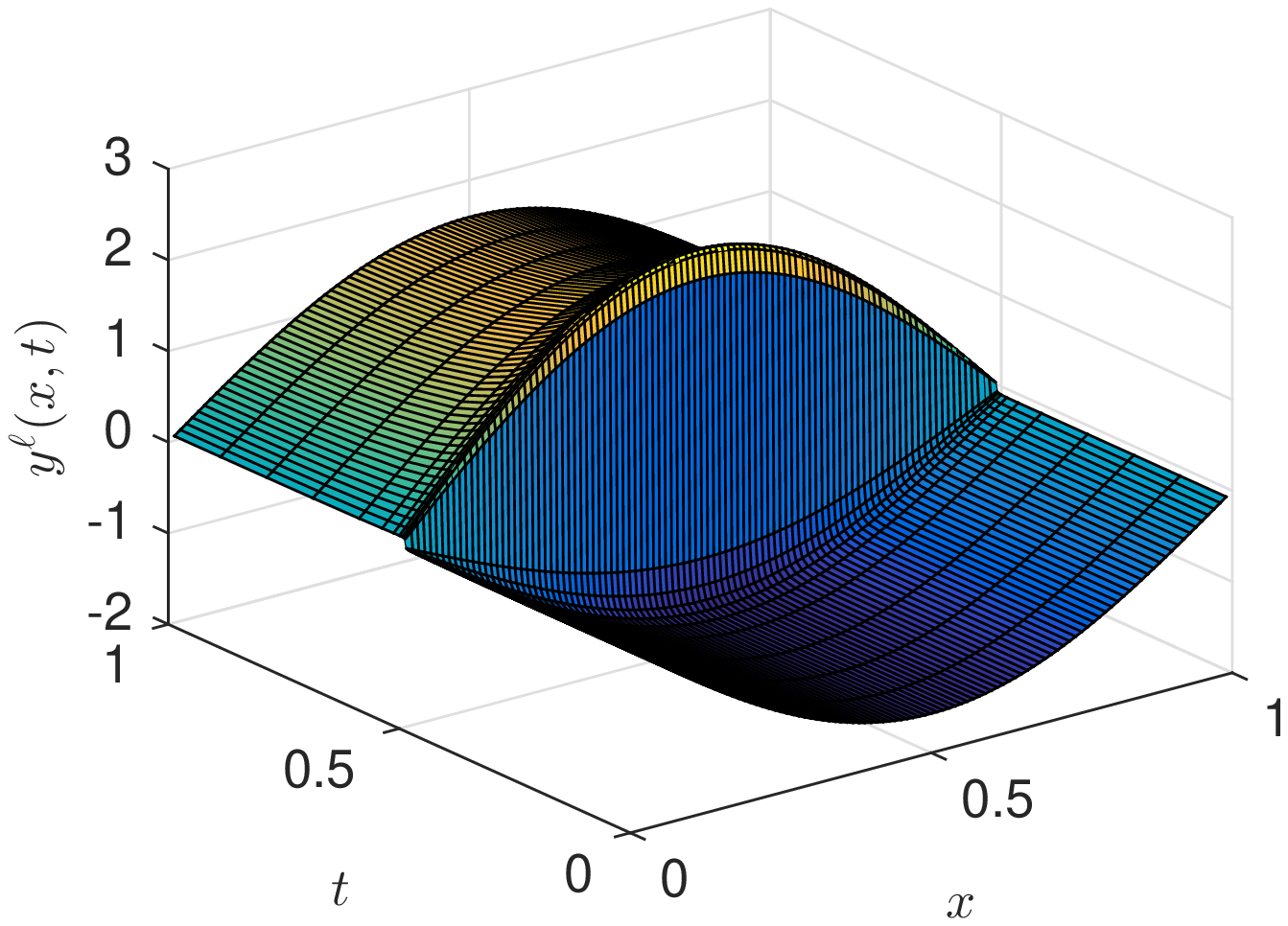}
 \caption{Exact optimal state $\bar{y}$ (top), approximated POD solution computed on an 
 equidistant time grid (middle), approximated POD solution utilizing the adaptive time grid (bottom)}
 \label{fig1:sol}
 \end{figure}

\begin{figure}[htbp]
 \centering
 \includegraphics[scale=0.5]{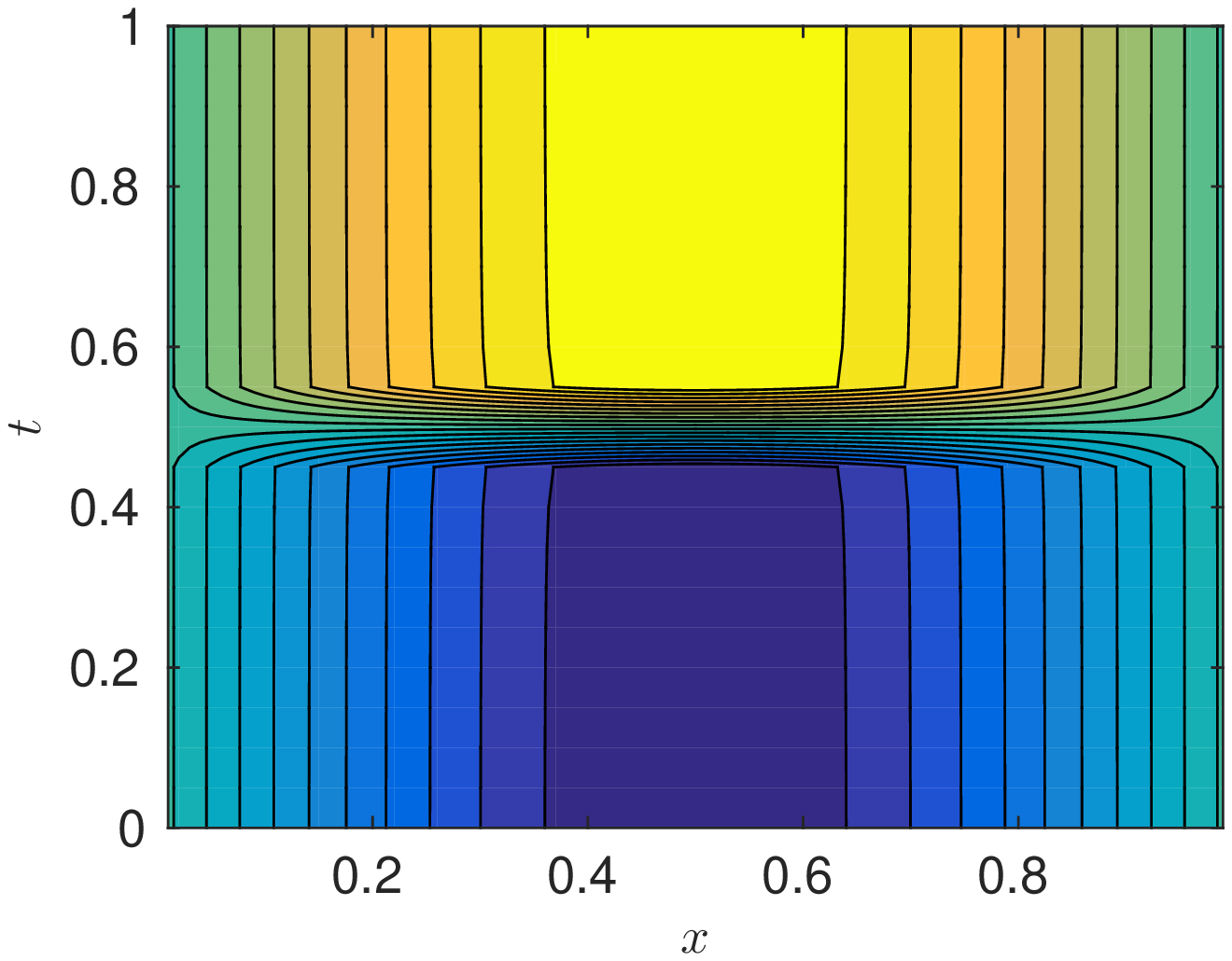}  \includegraphics[scale=0.5]{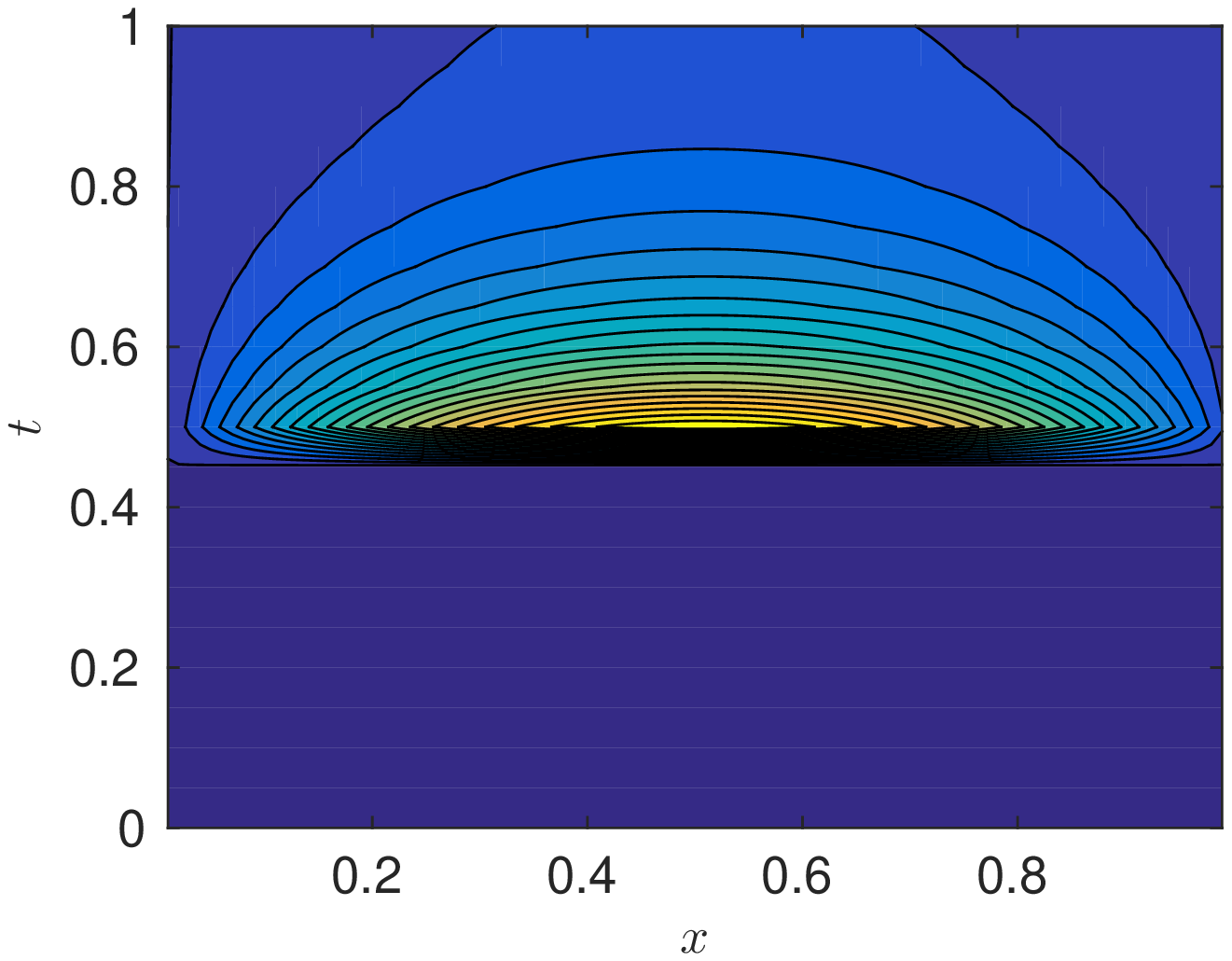}  \includegraphics[scale=0.5]{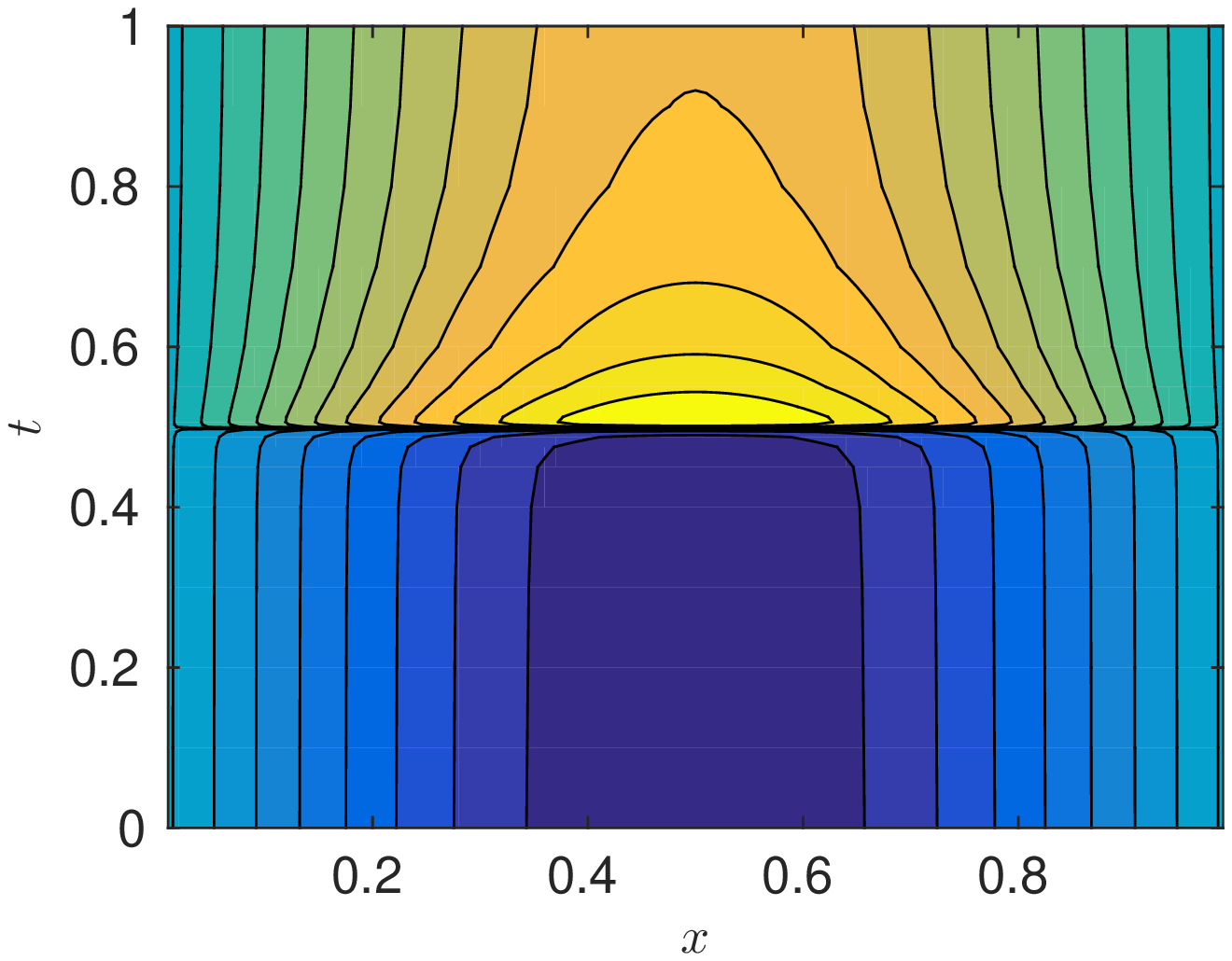}
 \caption{Contour lines of the exact optimal state $\bar{y}$ (top), 
 the approximated POD solution computed on an equidistant time grid (middle), the approximated POD 
 solution utilizing the adaptive time grid (bottom)}
 \label{fig1:con}
 \end{figure}

The initial condition is $y_0 (x) = \sin(\pi x) \text{atan}(-1/(2\varepsilon))$, $f$ and $y_d$ are chosen accordingly. For 
small values of $\varepsilon$ (we use $\varepsilon = 10^{-3}$), 
the state $\bar{y}$ develops an interior layer at $t = 1/2$, which can be 
seen at the top of Figure \ref{fig1:sol} and Figure \ref{fig1:con}.


\begin{figure}[htbp]
\centering
\includegraphics[scale=0.5]{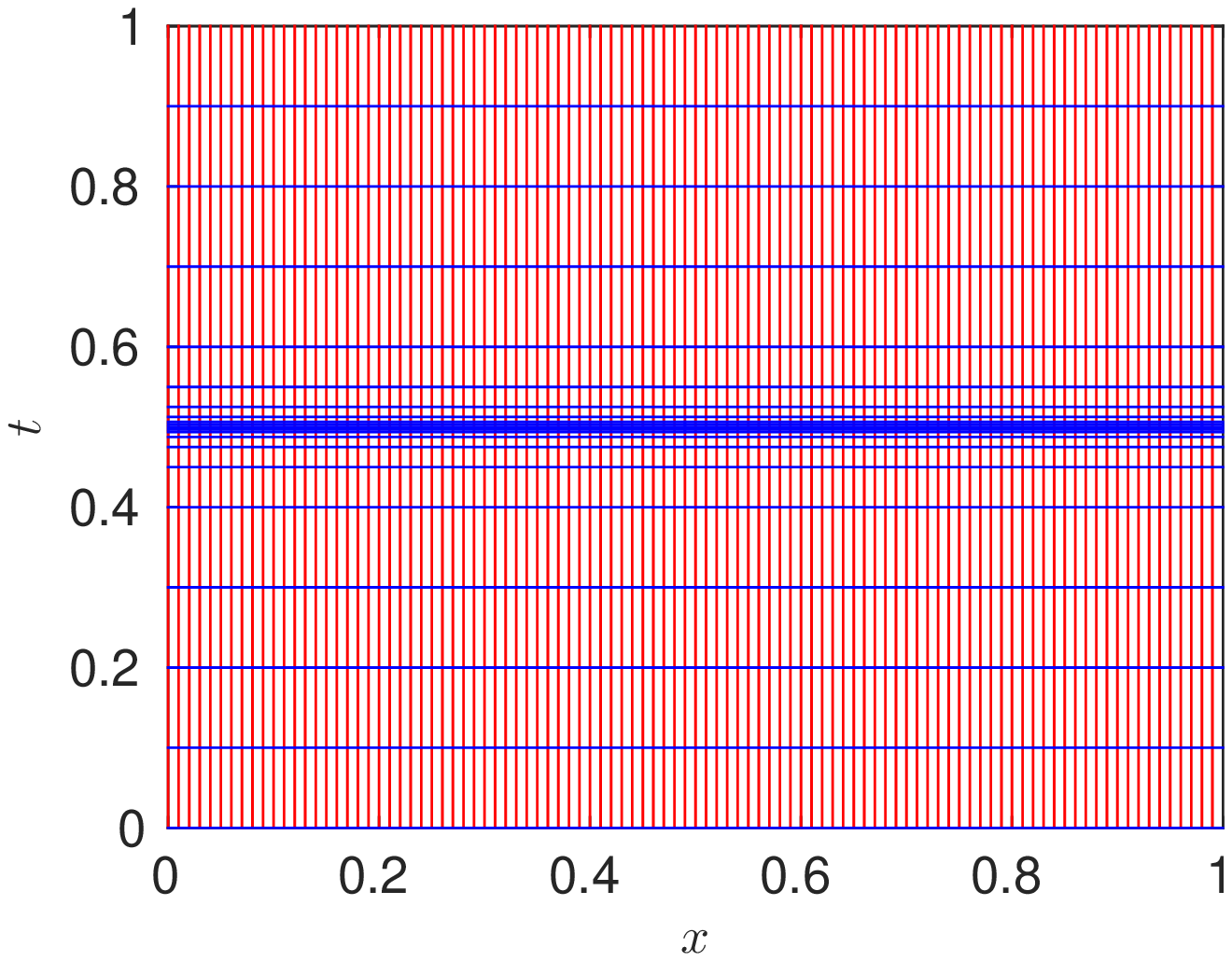}
\includegraphics[scale=0.5]{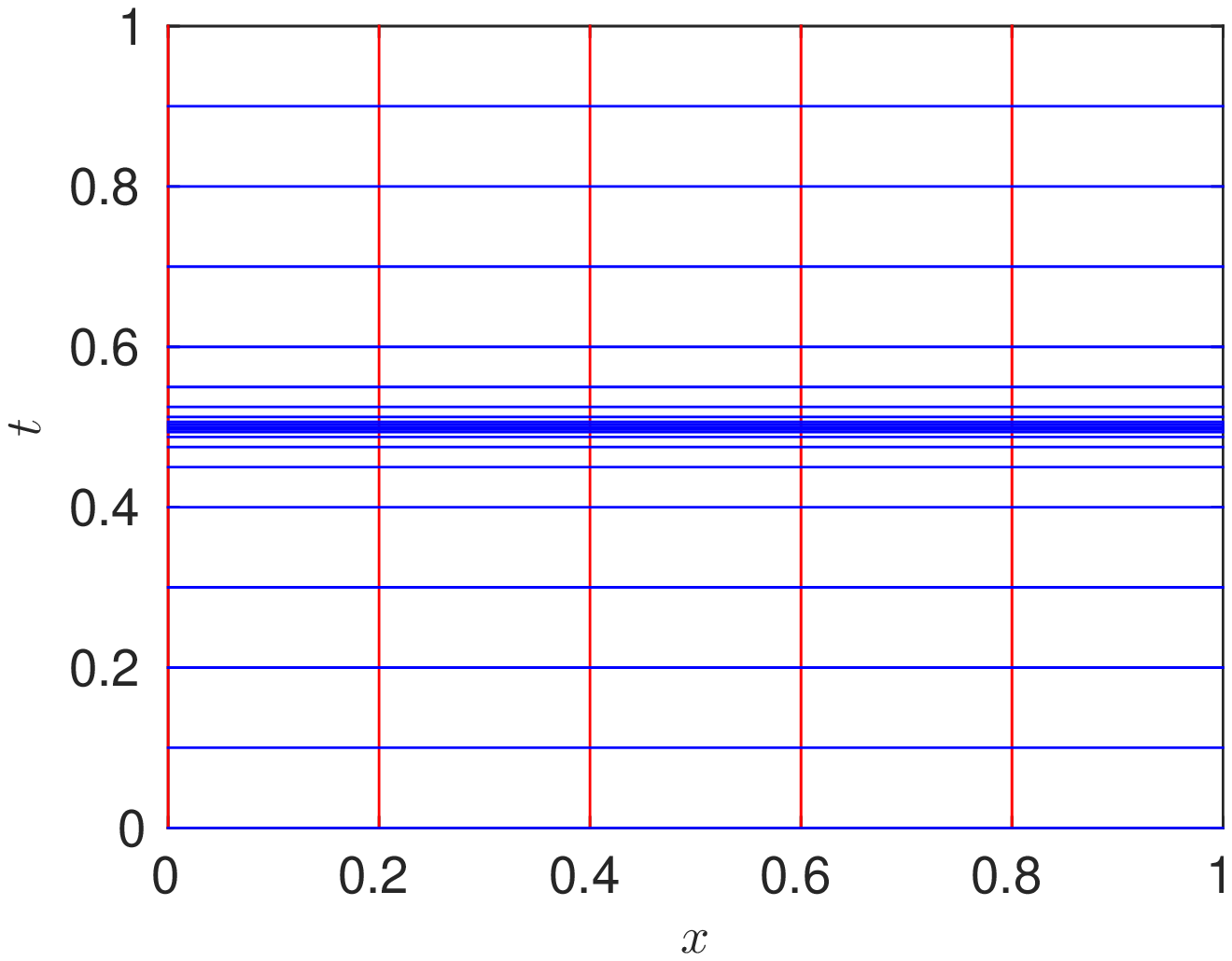}
\includegraphics[scale=0.5]{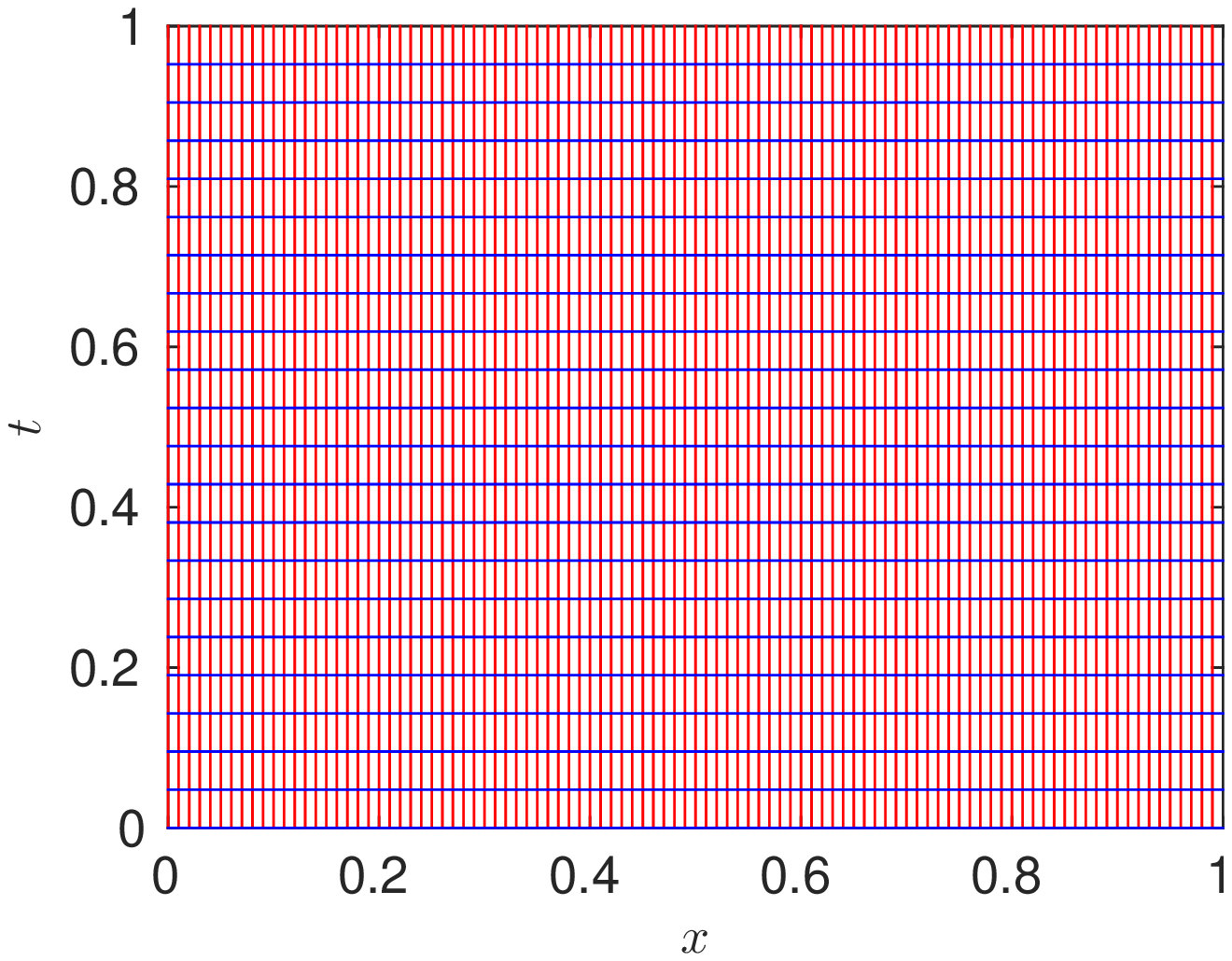}
    \caption{Adaptive space-time grids with dof $= 21$
    according to the strategy in \cite{GHZ12} and 
    $\Delta x = 1/100$ (top) and $\Delta x = 1/5$ (middle), 
    respectively, and the equidistant grid (bottom) with $\Delta x = 1/100$ and $\Delta t = 1/20$}
    \label{fig:1}
\end{figure}

 The main focus of our investigation consists in the use 
 of two different types of time 
grids: an equidistant time grid characterized by the 
time increment $\Delta t = 1/n$ and 
a non-equidistant (adaptive) time grid characterized 
by $n+1$ degrees of freedom (dof). Figure \ref{fig:1} visualizes the space-time mesh resulting from 
the strategy of \cite{GHZ12} utilizing temporal a-posteriori 
error estimation. The first grid in Figure \ref{fig:1} corresponds 
to the choice dof $= 21$ and $\Delta x = 1/100$, whereas the grid in the 
middle refers to using dof $= 21$ and $\Delta x = 1/5$. Both choices 
for spatial discretization lead to the exact same time grid, 
which displays fine time steps at time $t=1/2$ (where the layer in 
the optimal state is located), whereas at 
the beginning and end of the 
time interval the time steps are large. This clearly indicates that the resulting time-adaptive grid is very insensitive against changes in the spatial resolution. As a consequence we expect that the time grid obtained with a very coarse spatial resolution already delivers a time grid which well suits for the selection of fully resolved time-snapshots of the state, which form the basis of our POD-MOR strategy. For the sake of completeness, 
the equidistant grid with the same number of degrees of freedom is displayed at the bottom of Figure \ref{fig:1}.\\[1em]
\noindent Since the generation of the time-adaptive grid as well as the 
approximation of the optimal solution is done in the offline 
computation part of POD-MOR, this process shall be performed quickly, 
which is why we pick $\Delta x = 1/5$ for step 1 in Algorithm 
\ref{Alg:OPTPOD}.
For the computation of a POD suboptimal solution we make the 
following choices: we use uncontrolled snapshots, i.e. we 
create the snapshot ensemble by determining the 
associated state $y(u_{\text{sg}})$ to the control function 
$u_{\text{sg}} \equiv 0$. Although we would also 
have the possibility to use snapshots associated to a spatially coarse optimal 
control (since this optimal control is computed 
within Algorithm \ref{Alg:OPTPOD}), we want to emphasize here the 
importance of the time grid and compare with the most common snapshot selection approach in practice. 
Moreover, we use $1$ POD basis function determined 
by the \textsc{Matlab} routine $\mathtt{eigs}$, i.e. $\ell = 1$. 
As weighted inner product we choose 
$\langle v,w\rangle = \langle v, Ww \rangle_{\mathbb{R}^N}$ with $W=M$, where 
$N$ refers to the dimension of the finite element space associated with the spatial grid 
size $h$.\\
%
Table \ref{tab:1} and Table \ref{tab:1a} summarize all test 
runs with an equidistant and adaptive time grid, respectively. 
The fineness of time discretization is chosen 
in such a manner that the results are comparable. The absolute errors 
between the exact optimal state $\bar{y}$ and the 
POD suboptimal solution $\bar{y}^\ell$, defined by
$\varepsilon_{\text{abs}}^y := \parallel \bar{y} - \bar{y}^\ell
\parallel_{L^2(\Omega_T)}  $, are listed in columns 2,  same applies for the absolute errors in the control function, 
$\varepsilon_{\text{abs}}^u :=  \parallel \bar{u} - \bar{u}^\ell
\parallel_{L^2(\Omega_T)}$, 
columns 3.  If we compare the errors between exact  
and POD solution utilizing an equidistant time grid 
with the results for the errors utilizing the adaptive time grid, we clearly notice 
that the use of an adaptive time grid heavily improves the quality of the 
numerical POD solution.\\ 
\begin{table}[htbp]
\centering
 \begin{tabular}{ |c | c | c | }
 \toprule
 $\Delta t$ & $\varepsilon_{\text{abs}}^y$ 
 & $\varepsilon_{\text{abs}}^u$  \\
 \hline
 1/20 & $ 6.7264 \cdot 10^{+00}$  & $ 4.7581 \cdot 10^{-01}$  \\
 1/46 & $ 2.8634 \cdot 10^{+00}$ & $ 1.9863 \cdot 10^{-01}$ \\
 1/82 & $ 1.4603 \cdot 10^{+00}$ & $ 1.0037 \cdot 10^{-01}$ \\
 1/108 & $ 1.0233 \cdot 10^{+00}$ & $ 6.9951 \cdot 10^{-02}$ \\
 \bottomrule 
 \end{tabular}
 \vspace{0.4cm} \caption{Absolute errors between the exact optimal solution and 
 the POD suboptimal solution depending on the time 
 discretization with equidistant grid}
 \label{tab:1}
  \end{table}
  \begin{table}[htbp]
\centering
 \begin{tabular}{ | c | c | c | }
 \toprule
  dof & 
 $\varepsilon_{\text{abs}}^y$ & $\varepsilon_{\text{abs}}^u$ \\
 \hline
  21 & $ 1.8468 \cdot 10^{-01}$ & $ 1.1610 \cdot 10^{-02}$ \\
 47 & $ 2.0049 \cdot 10^{-02}$ & $ 9.3412 \cdot 10^{-03}$ \\
 83 & $ 1.2574 \cdot 10^{-02}$ & $ 5.0596 \cdot 10^{-03}$ \\
  109 & $ 6.3950 \cdot 10^{-03}$ & $ 4.8042 \cdot 10^{-03}$ \\
 \bottomrule 
 \end{tabular}
 \vspace{0.4cm} \caption{Absolute errors between the exact optimal solution and 
 the POD suboptimal solution depending on the time 
 discretization with adaptive grid }
 \label{tab:1a}
  \end{table}
 In fact, our numerical tests enable us to detect the 
 following: in order to achieve an accuracy in the state 
 variable of order 
 $10^{-1}$ utilizing an equidistant time grid, 
 we need about $n=110$ time steps and for an accuracy of 
 order $10^{-2}$ about $n=400$ time steps are 
 needed (not listed in Table 1). Once again, this emphasizes that using an 
 appropriate (non-equidistant) time grid is of particular 
 importance in order to efficiently achieve POD solutions of
 good quality.

 Figures \ref{fig1:sol} and \ref{fig1:con} (middle and bottom plots) show the surface and contour lines of the 
 POD suboptimal state utilizing an equidistant time grid (with 
 $\Delta t = 1/20$) and utilizing the adaptive 
 time grid (with \mbox{dof = 21}), respectively. As 
 expected, significant differences can be noticed in the appearance: 
 an equidistant time grid fails to capture 
 the interior layer at $t=1/2$ satisfactorily, whereas the POD solution utilizing the adaptive time grid 
 approxi\-mates the interior layer well.\\
 \noindent We note that enlarging the number of utilized POD basis 
 functions does not improve the approximation quality. 
Furthermore, we can argue that the  percentage of modelled energy to total energy contained in 
 the snapshots (energy content) is approximately 1.00 and the 
 second largest eigenvalue of the correlation matrix is of 
 order $10^{-16}$ (machine precision), which makes the use 
 of additional POD basis functions redundant. Likewise, as already mentioned, in 
 this particular example the choice of richer 
 snapshots (even the optimal snapshots) does not bring significant 
 improvements in the approximation quality of the POD 
 solutions. So, this example shows that 
 solely the use of an appropriate adaptive time mesh  
 effectively improves the accuracy of the POD suboptimal 
 solution. \\
 \noindent
The last remark goes to the efficiency of the POD method. 
 For the solution of the full-dimensional problem, 4 iterations 
 of the gradient method are needed, which means the 
 approximation of 8 (high-dimensional) PDEs is required, whereas the 
 reduced model took 5 iterations implying the approximation 
 of 10 ODEs of dimension $\ell = 1$. 
 The offline stage is really cheap since we use a very coarse spatial discretization.

 \section{Conlcusion and future directions}
In this paper we investigate the problem of snapshot location in optimal control problems. We show that the numerical POD solution is 
much more accurate if we use an adaptive time grid, especially when the solution 
of the problem presents layers. The time grid is computed by means of an a-posteriori error estimation 
strategy of the space-time approximation of a second order in time and fourth space elliptic equation which 
describes the optimal control problem and has the advantage that it does not explicitly depend on an input control function. 
Furthermore, a coarse approximation of the latter equation gives information on the snapshots one 
can use to build the surrogate model. 
Future work will focus on error analysis of the proposed algorithm, and on optimal control with shape functions and control constraints.



%
%

%


\begin{thebibliography}{4}
\providecommand{\natexlab}[1]{#1}
\providecommand{\url}[1]{\texttt{#1}}
\providecommand{\urlprefix}{URL }
\expandafter\ifx\csname urlstyle\endcsname\relax
  \providecommand{\doi}[1]{doi:\discretionary{}{}{}#1}\else
  \providecommand{\doi}{doi:\discretionary{}{}{}\begingroup
  \urlstyle{rm}\Url}\fi

\bibitem[{Able(1956)}]{Abl:56}
Able, B. (1956).
\newblock Nucleic acid content of microscope.
\newblock \emph{Nature}, 135, 7--9.

\bibitem[{Able et~al.(1954)Able, Tagg, and Rush}]{AbTaRu:54}
Able, B., Tagg, R., and Rush, M. (1954).
\newblock Enzyme-catalyzed cellular transanimations.
\newblock In A.~Round (ed.), \emph{Advances in Enzymology}, volume~2, 125--247.
  Academic Press, New York, 3rd edition.

\bibitem[{Keohane(1958)}]{Keo:58}
Keohane, R. (1958).
\newblock \emph{Power and Interdependence: World Politics in Transitions}.
\newblock Little, Brown \& Co., Boston.

\bibitem[{Powers(1985)}]{Pow:85}
Powers, T. (1985).
\newblock Is there a way out?
\newblock \emph{Harpers}, 35--47.

\end{thebibliography}


\begin{thebibliography}{xx} 
\bibitem[Afanasiev et~al.(2001)]{AH01}
K. Afanasiev and M. Hinze. 
\newblock Adaptive control of a wake flow using proper orthogonal decomposition.
\newblock \emph{Lecture Notes in Pure and Applied Mathematics 216}, 317-332. 
\newblock Shape Optimization \& Optimal Design, Marcel Dekker, 2001. 
\bibitem[Arian et~al.(2002)]{AFS02}
E. Arian, M. Fahl and E. Sachs.
\newblock Trust-region proper orthogonal decomposition models by optimization methods.
\newblock \emph{In Proceedings of the 41st IEEE Conference on Decision and Control, Las Vegas, Nevada}, pages 3300-3305, 2002.
\bibitem[Gong et~al.(2012)]{GHZ12}
W. Gong, M. Hinze, Z.J.Zhou.
\newblock Space-time finite element approximation of parabolic optimal control problems.
\newblock \emph{J. Numer. Math}, {\bf 20}, 2012, 111-145.
\bibitem[Gubisch et~al.(2013)]{GV13}
M. Gubisch and S. Volkwein.
\newblock Proper Orthogonal Decomposition for Linear-Quadratic Optimal Control.
\newblock \emph{Submitted}, 2013.
 \bibitem[Kunisch et~al.(2001)]{KV01} 
K. Kunisch, S. Volkwein.
\newblock Galerkin proper orthogonal decomposition methods for parabolic problems.
\newblock \emph{Numer. Math.} 90 (2001), 117-148.
\bibitem[Kunisch et~al.(2002)]{KV02}
K. Kunisch, S. Volkwein.
\newblock Galerkin proper orthogonal decomposition methods for a general equation in fluid dynamics.
\newblock \emph{SIAM, J. Numer. Anal.} {\bf 40} (2002), 492-515.
\bibitem[Kunisch et~al.(2008)]{KV08}
K. Kunisch and S. Volkwein.
\newblock Proper Orthogonal decomposition for optimality systems.
\newblock \emph{ESAIM: M2AN,} {\bf 42}, 2008, 1-23.
\bibitem[Kunisch et~al.(2010)]{KV10}
K. Kunisch and S. Volkwein.
\newblock Optimal Snapshot Location for computing POD basis functions.
\newblock \emph{ESAIM: M2AN,} {\bf 44}, 2010, 509-529.
\bibitem[Hinze (2005)]{H05}
M. Hinze.
\newblock A variational discretization concept in control constrained optimization: the linear-quadratic case.
\newblock \emph{Computational Optimization and Applications,} 30, 45-61, 2005
\bibitem[Hinze et~al.(2009)]{HPUU09} 
M. Hinze, R. Pinnau, M. Ulbrich and S. Ulbrich. 
\newblock Optimization with PDE Constraints. Mathematical Modelling: Theory and Applications, 23. 
\newblock \emph{Springer Verlag,} 2009.
\bibitem[Hoppe et~al.(2014)]{HL14}
R.H.W. Hoppe and Z. Liu.
\newblock Snapshot location by error equilibration in proper orthogonal decomposition for linear and semilinear parabolic partial differential equations
\newblock \emph{J. of Numer. Math.,} {\bf 22}, 2014, 1-32.

\bibitem[Neitzel et~al.(2009)]{NPS09}
I. Neitzel, U. Pr\"ufert, T. Slawig.
\newblock On Solving Parabolic Optimal Control Problems by Using Space-Time Discretization.
\newblock \emph{Technical Report} 05-2009, TU Berlin

\bibitem[Nocedal et~al.(2006)]{NW06} 
J. Nocedal and S.J. Wright. 
\newblock Numerical Optimization, second edition. 
\newblock \emph{Springer Series in Operation Research and Financial Engineering,} 2006.
%
\bibitem[Sirovich(1987)]{Sir87}
 L. Sirovich. 
\newblock Turbulence and the dynamics of coherent structures. Parts I-II.
\newblock \emph{Quarterly of Applied Mathematics,} {\bf XVL} (1987), 561-590.
\bibitem[Tr\"oltzsch(2010)]{Tro10}
F.~Tr\"oltzsch.
\newblock Optimal Control of Partial Differential Equations: Theory, Methods and Application.
\newblock \emph{American Mathematical Society,} 2010.
\bibitem[Tr\"oltzsch et~al.(2009)]{TV09}
F.~Tr\"oltzsch and S. Volkwein.
\newblock POD a-posteriori error estimates for linear-quadratic optimal control problems.
\newblock \emph{Comput. Optim. and Appl.,} {\bf 44}, 2009, 83-115.
%
\end{thebibliography}
\end{document}